\documentstyle[12pt]{article}
\begin{document}
\begin{center}
{\Large \bf Rigidity of surjective holomorphic maps }

\medskip
{\Large \bf to Calabi-Yau manifolds}

\bigskip
{\large \bf Jun-Muk Hwang} \footnote{Supported  by the Korea
Research Foundation Grant (KRF-2002-070-C00003).}
 \end{center}

\bigskip
\section{Introduction}

\medskip
Let $Y$ be a compact K\"ahler manifold with $c_1(Y) =0$ and $X$ be
a compact complex manifold. We denote by $\mbox{ Hol}_{surj}(X,
Y)$ the space of surjective holomorphic maps from $X$ to $Y$.
 In [KSW, Theorem 2], the following result  was proved.

 \medskip
 {\bf Theorem 0 (Kalka-Shiffman-Wong)} {\it Let $Y$ be a compact K\"ahler
 manifold of dimension $n$
 with $c_1(Y) =0$. Then for any compact complex manifold $X$,  $\mbox{Hol}_{surj}(X,Y)$
 is discrete if $c_n(Y) \neq 0$.}

 Although this result applies to many interesting examples such as K3-surfaces,
it is not completely satisfactory because the geometric nature of
K\"ahler manifolds satisfying $c_1(Y) =0$ and $ c_n(Y)=0$ is not
completely understood.

In this paper, we show that the idea of [KSW] can be easily
 generalized  to any $Y$ with $c_1(Y)=0$ by a slight refinement of the
argument. Of course, without any additional condition on $Y$,
$\mbox{Hol}_{surj}(X,Y)$ is not necessarily discrete. In fact, it
is possible that the identity component $A$ of the holomorphic
automorphism group of $Y$ has positive dimension, as in the case
of complex tori. Then the composition of a given element of
$\mbox{Hol}_{surj}(X,Y)$ with elements of $A$ gives a family of
surjective holomorphic maps. Our result says that, up to a finite
etale cover, this is the only way $\mbox{Hol}_{surj}(X, Y)$ fails
to be discrete. More precisely, we will prove the following.
Recall that $Y$ has a Ricci-flat K\"ahler metric by Yau's solution
of Calabi conjecture [Y].

\medskip
{\bf Theorem 1} {\it Let $Y$ be a compact Ricci-flat K\"ahler
manifold  and $X$ be a compact complex manifold. For any
surjective holomorphic map $f:X \rightarrow Y$, there exists a
finite etale cover $g:Z \rightarrow Y$ with the following
properties.

(1) There exists a surjective holomorphic map $h: X \rightarrow Z$
such that $f= g \circ h$.

(2) The two natural homomorphisms
$$H^0(Z, T(Z)) \stackrel{g_*}{\longrightarrow} H^0(Z, g^*T(Y))
\stackrel{h^*}{\longrightarrow} H^0(X, f^*T(Y))$$ are
isomorphisms. }

\medskip
Recall that the identity component of the holomorphic automorphism
group of a non-uniruled compact K\"ahler manifold is a complex
torus (e.g. [F, Theorem 5.5]). Since $H^0(X, f^*T(Y))$ is the
Zariski tangent space to $\mbox{Hol}_{surj}(X, Y)$ at the point
corresponding to $f$ (e.g., [H]), Theorem 1 implies the following.

\medskip
{\bf Corollary 1} {\it Let $f:X \rightarrow Y$ and $Z$ be as in
Theorem 1.
 Let $A$ be the complex torus which is the identity component of the holomorphic
automorphism group of $Z$. Then the holomorphic map $\psi: A
\rightarrow \mbox{\rm Hol}_{surj}(X,Y)$ defined by $$\psi(a) := g
\circ a \circ h$$ for $a \in A$ is unramified and surjective on
the component of $\mbox{\rm Hol}_{surj}(X, Y)$ containing $f$. In
particular, if $f$ is generically finite, holomorphic vector
fields on $X$ are pull-backs of those on $Z$. }

\medskip
 Thus  each component of $\mbox{Hol}_{surj}(X,Y)$ is either a single
point or an etale quotient of a complex torus. In particular,
infinitesimal deformations of $f$ are un-obstructed, even though
the obstruction group $H^1(X, f^*T(Y))$ may not be zero.

\medskip
Note that taking an etale cover is necessary in the statement of
Theorem 1. In fact, Igusa ([I, p.678]) constructed  an etale
quotient $f:X\rightarrow Y$ of a 3-dimensional abelian variety $X$
such that holomorphic vector fields on $X$ do not come from $Y$.

\medskip
Kalka, Shiffman and Wong noted (cf. the last paragraph of [KSW])
that Theorem 0 gives a partial generalization of Lichnerowicz'
result that the automorphism group of  a compact K\"ahler manifold
$Y$ with $c_1(Y)=0$  is discrete if the Chern numbers of $Y$ are
not all zero. Our result gives a full generalization:

\medskip
{\bf Corollary 2} {\it Let $Y$ be a compact K\"ahler manifold with
$c_1(Y) =0$. Then $\mbox{\rm Hol}_{surj}(X, Y)$ is discrete for
any compact complex manifold $X$, unless all the Chern numbers of
$Y$ vanish.}

\medskip
In fact, if a Chern number of $Y$ is non-zero, so is a Chern
number of the etale cover $Z$. Thus $Z$ cannot have a holomorphic
vector field as explained in the last paragraph of [KSW].

\medskip
It is well-known that a compact K\"ahler manifold $Y$ with
$c_1(Y)=0$ and $\pi_1(Y)$ finite has no holomorphic vector fields
(e.g., [B, p.759, Lemme]). Thus we have the following corollary.

\medskip {\bf Corollary 3} {\it Suppose  $Y$ is a compact K\"ahler
 manifold with $c_1(Y) =0$ and $|\pi_1(Y)|< +\infty$, e.g. a Calabi-Yau manifold.
 Then $\mbox{\rm Hol}_{surj}(X,Y)$ is discrete for any compact complex manifold
 $X$.   }

\medskip
After the first version of this paper was written, the following
generalization of  Theorem 1 was proved in [HKP] when $Y$ is an
algebraic variety.

\medskip
{\bf Theorem [HKP]}  {\it Let $Y$ be a projective algebraic
manifold with non-negative Kodaira dimension, or more generally, a
projective algebraic manifold which is not covered by rational
curves. Let $X$ be a compact complex manifold. For any surjective
holomorphic map $f:X \rightarrow Y$, there exists a finite etale
cover $g:Z \rightarrow Y$ with the following properties.

(1) There exists a surjective holomorphic map $h: X \rightarrow Z$
such that $f= g \circ h$.

(2) The two natural homomorphisms
$$H^0(Z, T(Z)) \stackrel{g_*}{\longrightarrow} H^0(Z, g^*T(Y))
\stackrel{h^*}{\longrightarrow} H^0(X, f^*T(Y))$$ are
isomorphisms. }

\medskip
Corresponding generalizations of Corollary 1 and Corollary 3
follow.  The method used in [HKP] is completely different from the
current paper and depends heavily on the assumption that $Y$ is
algebraic. Whether the same statement holds for all compact
K\"ahler manifolds of non-negative Kodaira dimension, which would
be a full generalization of Theorem 1, remains open.

\section{ Proof of Theorem 1}

\medskip
 The proof of Theorem 1 is based on the following two lemmata.
 Lemma 1 is contained in the proof of [KSW, Theorem 2], but we will give a full proof
 because the notation and an inequality there  will be used in  the
 proof of Lemma 2.

\medskip {\bf Lemma 1} {\it Let $j: W \rightarrow Y$ be a finite
holomorphic map from a complex variety $W$ to a compact complex
manifold $Y$ with a Ricci-flat K\"ahler metric $\omega$. Then for
any section $\sigma \in H^0(W, j^* T(Y))$, its norm $||\sigma||$
with respect to the metric on $j^*T(Y)$ induced by $\omega$ is a
constant function on $W$.}

\medskip
{\it Proof}. For a given non-zero section $\sigma \in H^0(W,
j^*T(Y))$, let $a = \sup_{W} ||\sigma||
>0$ and let $w_o \in W$ such that $||\sigma_{w_o}|| = a$. It
suffices to show that $||\sigma|| \equiv a$ in a neighborhood of
$w_o$ in $W$. Let $y_o = j(w_o).$ Choose a connected neighborhoods
$W_o$ of $w_o$ and $Y_o$ of $y_o$ such that $$j_o:= j|_{W_o}: W_o
\rightarrow Y_o$$ is a $\lambda$-sheeted branched analytic cover
for some $\lambda \geq 1$ with $j_o^{-1}(y_o) = w_o$
set-theoretically. Define a holomorphic vector field
$\bar{\sigma}$ on $Y_o$ by setting its value at $y \in Y_o$ to be
$$\bar{\sigma}_y := \sum_{i=1}^{\lambda} \sigma_{w^i}$$ where
$w^1, \ldots, w^{\lambda}$ are the points of $j^{-1}(y)$ counted
with multiplicities and $\sigma_{w^i}$ denotes the value of
$\sigma$ at $w^i$. In particular,
$$\bar{\sigma}_{y_o} = \lambda \sigma_{w_o}.$$ For any $y \in
Y_o$,
$$||\bar{\sigma}_y|| \leq \sum_{i=1}^{\lambda} ||\sigma_{w^i}||
\leq \lambda a = ||\bar{\sigma}_{y_o}||.$$ Thus
$||\bar{\sigma}||^2$ attains its maximum at $y_o$. Now we recall
the following.

\medskip
{\bf Bochner formula ([B, p.760] or [KSW, Lemma 4])} {\it Let $s$
be a holomorphic tensor on a Ricci-flat K\"ahler manifold. Then
$\Delta (||s||^2) = ||\nabla s||^2.$}

\medskip
Applying it to $s= \bar{\sigma}$, we have
$$\Delta (||\bar{\sigma}||^2) = ||\nabla \bar{\sigma}||^2 \geq 0$$
on $Y_o$. Hence $||\bar{\sigma}||^2$ is subharmonic on $Y_o$.
Since it attains its maximum at an interior point $y_o$,
$||\bar{\sigma}||^2$ must be constant on $Y_o$. Since for $y \in
Y_o$, $$(\dagger) \;\;\;\;\; \lambda a \geq \sum_{i=1}^{\lambda}
||\sigma_{w^i}|| \geq ||\bar{\sigma}_y|| = \lambda a,$$ it follows
that $||\sigma|| \equiv a $ on $W_o$. $\Box$

\medskip
{\bf Lemma 2} {\it Let $j: W_o \rightarrow Y_o$ be a finite
holomorphic map of degree $\lambda$ from an irreducible complex
variety $W_o$ to a complex manifold $Y_o$ having a Ricci-flat
K\"ahler metric $\omega$. Suppose for a point $y_o \in Y_o$, the
set-theoretic inverse image $j^{-1}(y_o)$ is a single point $w_o$.
Let $\sigma \in H^0(W_o, j^* T(Y_o))$ be a section such that its
norm $||\sigma||$ with respect to the metric on $j^*T(Y_o)$
induced by $\omega$ is a constant function on $W_o$. Then for any
point $y \in Y_o$ and $j^{-1}(y) = \{ w^1, \ldots, w^{\lambda}\}$
counting multiplicities,
$$\sigma_{w^1} = \sigma_{w^2} = \cdots = \sigma_{w^{\lambda}}$$
regarded as vectors in $T_y(Y_o)$.}

\medskip
{\it Proof}.  Since $||\sigma||$ is constant, we may regard
$||\sigma_{w_o}||$ as the maximum value of $||\sigma||$ and
 apply the argument in the proof of Lemma 1.  By the equality
$(\dagger)$ in the proof of Lemma 1,
$$\sum_{i=1}^{\lambda} ||\sigma_{w^i}|| = \lambda
||\sigma_{w_o}||$$ for any $y \in Y_o$ and $j^{-1}(y) = \{ w^1,
\ldots, w^{\lambda} \}$. This implies the desired equality
$$\sigma_{w^1} = \sigma_{w^2} = \cdots =
\sigma_{w^{\lambda}}$$ by the following elementary fact. $\Box$

\medskip
{\bf Sublemma} {\it
 Let $v^1, \ldots, v^{\lambda}$ be elements of a complex
vector space equipped  with  a Hermitian inner product and
corresponding norm $|| \cdot ||$. If $||v^1|| = ||v^2|| = \cdots =
||v^{\lambda}|| $ and $|| \sum_{i=1}^{\lambda} v^{i} || = \lambda
||v^1||,$ then $v^1 = v^2 = \cdots = v^{\lambda}$. }

\medskip
{\it Proof}.  This follows from the  identity
$$ ||\sum_{i=1}^{\lambda} v^i||^2 + \sum_{i<k}||v^i-v^k||^2
= \lambda \sum_{i=1}^{\lambda} ||v^i||^2. \;\;\; \Box$$

\medskip
Now let $X$ be a compact complex manifold and  $f: X \rightarrow
Y$ be a surjective holomorphic map. Let $\hat{f}: \hat{X}
\rightarrow Y$ and $\tilde{f}: X \rightarrow \hat{X}$ be the Stein
factorization of $f$. Then we have a natural isomorphism
$$H^0(X, f^*T(Y)) \cong H^0(\hat{X}, \hat{f}^* T(Y))$$ because $\tilde{f}$ has connected fibers.
 We define
an equivalence relation $R$ on $\hat{X}$ as follows. We say that
two points $x$ and $x'$ of $\hat{X}$ are $R$-equivalent if

(a) $\hat{f}(x) = \hat{f}(x')$ and

(b) for each $\sigma \in H^0(\hat{X}, \hat{f}^*T(Y))$, $\sigma_x =
\sigma_{x'}$ as vectors in $T_{\hat{f}(x)}(Y) =
T_{\hat{f}(x')}(Y)$.

\medskip We claim that the quotient $\hat{X} \rightarrow \hat{X}/R$ of $\hat{X}$ by the
equivalence relation $R$ exists as a holomorphic map to a complex
analytic variety. Denote by $$\hat{f}_*: \hat{f}^*T(Y) \rightarrow
T(Y)$$ the finite holomorphic map induced by $\hat{f}$ between the
total spaces of vector bundles. Viewing each $\sigma \in
H^0(\hat{X}, \hat{f}^*T(Y))$ as a holomorphic map
$$\sigma: \hat{X} \rightarrow \hat{f}^*T(Y),$$  define
$X_{\sigma} \subset T(Y)$ by $$X_{\sigma} :=
\hat{f}_*(\sigma(\hat{X})).$$ The natural holomorphic map
$h_{\sigma}: \hat{X} \rightarrow X_{\sigma}$ defined by the
composition of $\sigma$ and $\hat{f}_*$ identifies two points $x$
and $x'$ of $\hat{X}$ satisfying (a) and $\sigma_x = \sigma_{x'}$
in the sense of (b). Thus the equivalence relation $R$ is defined
by the family of holomorphic maps $$\{ h_{\sigma}: \hat{X}
\rightarrow X_{\sigma} \; | \; \sigma \in H^0(\hat{X},
\hat{f}^*T(Y)) \}$$ in the sense of [C,p.8, Remark]. Then the
quotient $\hat{X} \rightarrow \hat{X}/R$ exists as a complex
analytic variety by [C,p.7, Main Theorem]. This proves the claim.

\medskip
Let $Z$ be the normalization of  $\hat{X}/R$. By (a), there exists
a natural finite holomorphic map $g: Z \rightarrow Y$. By (b),
there exists a natural isomorphism
$$H^0(X, f^*T(Y)) \cong H^0(\hat{X}, \hat{f}^*T(Y)) \cong H^0(Z,
g^*T(Y)).$$ By (b) and Lemma 2, $g$ must be an etale morphism. It
follows that $Z$ is smooth and $T(Z) = g^* T(Y)$, which implies
the desired equality
$$H^0(X, f^*T(Y)) \cong H^0(Z, g^*T(Y)) = H^0(Z, T(Z)).$$

\bigskip
{\bf Acknowledgment} I would like to thank N. Mok, B. Shiffman,
W.-K. To and B. Wong for helpful discussions and  encouragement.

\bigskip
 {\bf References}

\medskip
 [B] Beauville, A.: Vari\'et\'es K\"ahleriennes dont la
premi\`ere classe de Chern est nulle. J. Diff. Geom. {\bf 18}
(1983) 755-782

[C] Cartan, H.: Quotients of complex analytic spaces. {\it
International Colloquium on Function Theory.}  Tata Institute
(1960) 1-15 (equivalently, Article 51 in Collected Works Volume 2,
Springer Verlag, 1979)

 [F] Fujiki, A.: Automorphism groups of compact K\"ahler
manifolds. Invent. math. {\bf 44} (1978) 225-258

[H] Horikawa, E.: On deformations of holomorphic maps I. J. Math.
Soc. Japan, {\bf 25} (1973) 372-396

[HKP] Hwang, J.-M., Kebekus, S. and Peternell, T.: Holomorphic
maps onto varieties of non-negative Kodaira dimension. preprint.

 [I] Igusa, J: On the structure of a certain class of
K\"ahler varieties. Amer. J. Math. {\bf 76} (1954) 669-678

 [KSW] Kalka, M., Shiffman, B. and Wong, B.: Finiteness and
rigidity theorems for holomorphic mappings. Mich. Math. J. {\bf
28} (1981) 289-295

[Y] Yau, S.T.: On the Ricci curvature of a compact K\"ahler
manifold and the complex Monge-Amp\`ere equation. I. Comm. Pure
Appl. Math. {\bf 31} (1978) 339-411
\bigskip

\bigskip
Korea Institute for Advanced Study

207-43 Cheongnyangni-dong

Seoul, 130-722,  Korea

jmhwang@kias.re.kr
\end{document}